\documentclass[11pt,twoside]{amsart}

\usepackage{amssymb}
\usepackage{latexsym}

\newcommand{\Id}{\mathrm{Id}}

\newcommand{\dopu}{{:}\allowbreak\ }
\newcommand{\eps}{\varepsilon}
\newcommand{\cal}{\mathcal}

\newcommand{\mytilde}[1]{\mathbin{\tilde#1}}

\newcommand{\loglike}[1]{\mathop{\rm #1}\nolimits}
\newcommand{\ex}{\loglike{ex}}

\newcommand{\dist}{\loglike{dist}}

\newcommand{\lin}{\loglike{lin}}
\newcommand{\conv}{\loglike{conv}}

\newcommand{\Daug}{\loglike{Daug}}

\newcommand{\SDt}{\loglike{{\cal S}\!{\cal D}}}
\newcommand{\narr}{\loglike{{\cal N}\!\!{\cal A}\!{\cal R}}}

\newcommand{\Nar}{\narr}
\newcommand{\SD}{\SDt}
\newcommand{\rNar}{\Nar^r}
\newcommand{\rSD}{\SD^r}
\def\DP{Daugavet property}
\def\DPr{\mathrm{DP}}
\def\rDPr{\mathrm{DP}^r}
\newcommand{\sDo}{strong Daugavet operator}

\newcommand{\N}{{\mathbb N}}
\newcommand{\R}{{\mathbb R}}

\newcommand{\calU}{\mathcal{U}}

\theoremstyle{plain}
\newtheorem{thm}{Theorem}[section]

\newtheorem{prop}[thm]{Proposition}
\newtheorem{cor}[thm]{Corollary}
\newtheorem{lemma}[thm]{Lemma}

\theoremstyle{definition}
\newtheorem{definition}[thm]{Definition}

\theoremstyle{remark}

\newtheorem{rem}[thm]{Remark}

\numberwithin{equation}{section}

\newcommand{\begsta}{\begin{statements}}
\newcommand{\begaeq}{\begin{aequivalenz}}
\def\endsta{\end{statements}}
\def\endaeq{\end{aequivalenz}}
\newcommand{\bea}{\begin{eqnarray*}}
\newcommand{\eea}{\end{eqnarray*}}
\newcommand{\kref}[1]{(\ref{#1})}

\newcounter{abc}   
\newcounter{iiiii} 

\newenvironment{aequivalenz}
{\setcounter{iiiii}{0}
\begin{list}%
{{\rm (\roman{iiiii})}}
{\usecounter{iiiii}
\parsep=0pt plus 1pt
\topsep=1pt plus 2pt minus 1pt
\itemsep=1pt plus 2pt minus 1pt
\leftmargin=3\baselineskip
\labelsep=.6\baselineskip
\labelwidth=2.4\baselineskip
\rightmargin 0pt}%
}%
{\end{list}}

\newenvironment{statements}%
{\setcounter{abc}{0}
\begin{list}%
{{\rm (\alph{abc})}}
{\usecounter{abc}
\parsep=0pt plus 1pt
\topsep=1pt plus 2pt minus 1pt
\itemsep=1pt plus 2pt minus 1pt
\leftmargin=3\baselineskip
\labelsep=.6\baselineskip
\labelwidth=2.4\baselineskip
\rightmargin 0pt}%
}%
{\end{list}}

\begin{document}

\title [The Daugavet property for ultraproducts ]{ Narrow 
         operators and the Daugavet property for ultraproducts } 

\author[D.~Bilik, V.~Kadets, R.~Shvidkoy, D.~Werner]
{ Dmitriy Bilik, Vladimir Kadets, Roman Shvidkoy, \\
  and Dirk Werner }

\address{Faculty of Mechanics and Mathematics, Kharkov National
University,\linebreak
 pl.~Svobody~4,  61077~Kharkov, Ukraine}

\address{Faculty of Mechanics and Mathematics, Kharkov National
University,\linebreak
 pl.~Svobody~4,  61077~Kharkov, Ukraine}
\email{vishnyakova@ilt.kharkov.ua}

\address{Department of Mathematics, University of Missouri,
Columbia MO 65211}
\email{shvidkoy\_r@yahoo.com}

\address{Department of Mathematics, Freie Universit\"at Berlin,
Arnimallee~2--6, \qquad {}\linebreak D-14\,195~Berlin, Germany}
\email{werner@math.fu-berlin.de}

\thanks{The work of the second-named author
was supported by a grant from the {\it Alexander-von-Humboldt
Stiftung}.}




\begin{abstract}
We show that if $T$ is a narrow operator (for the definition see
below) on $X=X_{1}\oplus_{1} X_{2}$ or $X=X_{1}\oplus_{\infty}
X_{2}$, then the restrictions to $X_{1}$ and $X_{2}$ are narrow and
conversely. We also characterise by a version of the Daugavet property
for positive operators on Banach lattices which unconditional sums of Banach
spaces inherit the Daugavet property, and we study the Daugavet
property for ultraproducts.
\end{abstract}

\maketitle


\section{Introduction}
\label{sec1}

A Banach space $X$ is said to have the \textit{Daugavet
property} if every
rank-$1$ operator $T\dopu X \to X$ satisfies
\begin{equation}\label{eq1.1}
\|\Id+T\| =1+ \|T\|.
\end{equation}
Examples include $C[0,1]$, $L_{1}[0,1]$, certain function algebras
such as the disk algebra or $H^\infty$ and also their noncommutative
counterparts (nonatomic $C^*$-algebras and preduals of
nonatomic von Neumann algebras 
\cite{Oik}). Such spaces are studied in detail in \cite{KadSSW}.

It has long been known that the $\ell_{1}$-sum and the
$\ell_{\infty}$-sum of two spaces with the \DP\ again enjoys the \DP;
see \cite{Abra2}, \cite{Woj92} and \cite{KadSSW} for various degrees
of generality in this statement. In this paper we present an
operator theoretic extension of this result along with a study of the
inverse problem of which unconditional direct sums inherit the \DP.
Our methods are naturally related to ultrapower techniques, so we
also study the \DP\ for ultraproducts of Banach spaces. 

In \cite{KadSW2} we have introduced the notions of a \sDo\ and of a
narrow operator between Banach spaces $X$ and $Y$; 
for definitions see Section~\ref{sec2}. These appear to be,
in the case $X=Y$, the largest reasonable classes of operators that
satisfy \kref{eq1.1}; in particular, weakly compact operators,
operators not fixing a copy of $\ell_{1}$, strong Radon-Nikod\'ym
operators and their sums are \sDo s on Banach spaces with the \DP. 
In Section~\ref{sec3} we show
that an operator $T$ on $X_{1}\oplus_{\infty} X_{2}$ is a \sDo\ if
and only if both restrictions of $T$ to $X_{1}$ resp.\ $X_{2}$ are
\sDo s; the same is true for narrow operators. Section~\ref{sec4}
studies the same problem on $\ell_{1}$-sums $X_{1}\oplus_{1} X_{2}$;
here a subtle difference between \sDo s and narrow operators is
exhibited, since the restriction of a narrow operator turns out to be
narrow again, but the analogue for \sDo s proves to be false. 

Section~\ref{sec5} deals with the inverse problem to determine which
$1$-uncondi\-tional sums of Banach spaces inherit the
\DP. These are completely classified, and it turns out that among all
$1$-uncondi\-tional sums of two spaces only $X_{1} \oplus_{ 1} X_{2}$
and $X_{1} \oplus_{ \infty} X_{2}$ are spaces with the \DP.

As pointed out above, our methods rely on ultrapower techniques as
explained in Section~\ref{sec2}. Thus it appears natural to try and
investigate which ultrapowers $X^\calU$  have the \DP. These can be
characterised by means of a quantitative version of the \DP\ for $X$ that we
call the uniform \DP. As examples we show that $C[0,1]$ and
$L_{1}[0,1]$ have the uniform \DP.

We consider real Banach spaces in this paper; $S(X)$ stands for the
unit sphere of a Banach space $X$ and $B(X)$ for its closed unit ball.


\section{The rigid versions of the Daugavet property, strong 
Daugavet and narrow operators}
\label{sec2}

We recall the following characterisation of the Daugavet
property from \cite{KadSSW}:

\begin{lemma} \label{l:MAIN}
The following assertions are equivalent:
\begin{aequivalenz}
\item
     $X$ has the Daugavet property.
\item
     For every $x\in S(X)$, $x^*\in S(X^*)$  and $\eps >0$ there
         exists some $y\in S(X)$  such that $x^*(y) \ge 1-\eps$  and
$\|x+y \| \ge 2-\eps$.
\end{aequivalenz}
\end{lemma}

An a bit more general notion will be convenient for us.

\begin{definition}
A Banach space $X$ has the Daugavet
property with respect to a subset  $\Gamma \subset S(X^*)$ 
($X\in \DPr(\Gamma)$ for short) if for every $x\in S(X)$, $x^*\in \Gamma$  
and $\eps >0$ there  exists some $y\in S(X)$  such that $x^*(y)\ge1-\eps$  and
$\|x+y\| \ge 2-\eps$.
\end{definition}

According to \cite{KadSW2} 
an operator $T$ on a Banach space $X$ is said to be a {\em strong
Daugavet operator} if for every two elements $x, y \in S(X)$ and
for every $\eps > 0$ there is an element $z \in S(X) $ such that
$\|x+z\| \ge 2 - \eps$ and $\|T(y-z)\| \le \eps$. We denote the set of
all strong Daugavet operators on $X$ by $\SD(X)$.
An  operator $T$  is said to be a {\em narrow
operator} if for every two elements $x, y \in S(X)$, for every
$x^* \in X^*$ and for every $\eps > 0$ there is an element $z \in
S(X) $ such that $\|x+z\| \ge 2 - \eps$ and
$\|T(y-z)\|+|x^*(y-z)| \le \eps$. We denote the set of all narrow
operators on $X$ by $\Nar(X)$.

To indicate the difference between the two classes we have introduced
the following notation in \cite{KadSW2}. If $T\dopu X\to Y$ is an
operator and $x^*\dopu X\to \R$ is a functional, define
$$
T \mytilde+ x^*\dopu X\to Y\oplus _{1} \R, \quad x\mapsto (Tx,
x^*(x)).
$$
Then $T$ is narrow if and only if $T\mytilde+ x^*$ is a \sDo\ for
every $x^*\in X^*$.

A $\Gamma$-version of the definition of a narrow operator
will be also useful for us.

\begin{definition} \label{narG}
An  operator $T$ on a Banach space $X$ is said to be narrow with 
respect to a subset  $\Gamma \subset S(X^*)$
 ($T \in \Nar(X,\Gamma)$ for short) if for every two elements 
$x, y \in S(X)$, for every
$x^* \in \Gamma$ and for every $\eps > 0$ there is an element $z \in
S(X) $ such that $\|x+z\|\ge 2 - \eps$ and
$\|T(y-z)\|+|x^*(y-z)| \le \eps$. 
\end{definition}

It will be technically convenient to work with the case of $\eps=0$
in the above definitions. Therefore
we introduce ``rigid versions'' of these notions.

\begin{definition} \label{rigiddef}
\mbox{}
\begsta
\item 
A Banach space $X$ has the rigid Daugavet
property with respect to a subset  $\Gamma \subset S(X^*)$ 
($X\in \rDPr(\Gamma)$ for short) if for every $x\in S(X)$ and 
$x^*\in \Gamma$  
there  exists some $y\in S(X)$  such that $x^*(y)=1$  and
$\|x+y\|=2$.
\item
An operator $T$ on a Banach space $X$ is said to be a rigid strong
Daugavet operator (in symbols $T \in \rSD(X)$) if for every two 
elements $x, y \in S(X)$ there is an element $z \in S(X) $ such that
$\|x+z\|=2$ and $T(y-z)=0$.
\item
An  operator $T$  is said to be rigidly narrow with 
respect to a subset  $\Gamma \subset S(X^*)$
 (in symbols $T \in \rNar(X,\Gamma)$) if for every two elements 
$x, y \in S(X)$ and  for every
$x^* \in \Gamma$  there is an element $z \in
S(X) $ such that $\|x+z\|= 2$ and
$\|T(y-z)\|+|x^*(y-z)|=0$. 
\endsta
\end{definition}

Let us mention that  a rigid strong
Daugavet operator is necessarily non-injective. To see this, 
one just has to apply the 
definition for $y=-x$; then $y-z$ will be a nonzero element which
$T$ maps to 0. Using this remark one can easily prove the following
statement.

\begin{lemma} \label{equiv}
If $T \in \rSD(X)$, then  for every $x  \in S(X)$ and  
$y \in B(X)$ there is an element  $z \in S(X) $ such that
$\|x+z\|=2$ and $T(y-z)=0$.
\end{lemma}

\begin{proof}
Using the non-injectivity of $T$ one can find an element 
$y_1 \in S(X)$ such that $T(y-y_1)=0$. Then applying the definition 
of $\rSD(X)$ to $x$ and $y_1$ one obtains an element  
$z \in S(X) $ such that
$\|x+z\|=2$ and $T(y_1-z)=0$. But for this element $T(y-z)=0$, too.
\end{proof}

For many investigations in the context of the \DP\ the study of the
rigid notions above turns out to be sufficient, but is technically
more feasible. The connection between the original versions and their
rigid variants is made using  ultrapowers. We refer for
instance to \cite{Beau} for an introduction to ultrapowers of Banach
spaces.

Let $\cal U$ be a nontrivial ultrafilter on $\N$, $T$ be an operator
acting from a Banach space $X$ to a Banach space $Y$,  
$\Gamma \subset S(X^*)$. We denote by
$T^{\cal U}$ the natural operator between the ultrapowers $X^{\cal U}$
and $Y^{\cal U}$ defined by $T^{\cal U}(x_n)=(Tx_n)$, and by $\Gamma^{\cal U}$
we denote the set of the linear functionals 
$F=(f_n)$, $f_n \in \Gamma$, of the form $F(x_n)= \lim_{\cal U}f_n(x_n)$.

\begin{lemma} \label{rigid}
\mbox{ }
\begin{enumerate}
\item \label{r1}
If $X\in \DPr(\Gamma)$, then $X^{\cal U} \in \rDPr(\Gamma^{\cal U})$.
\item \label{r2}
If $X^{\cal U} \in \DPr(\Gamma^{\cal U})$, then $X\in \DPr(\Gamma)$.
\item \label{r3}
If $T \in \SD(X)$, then $T^{\cal U} \in \rSD(X^{\cal U})$.
\item \label{r4}
If $T^{\cal U} \in \SD(X^{\cal U})$, then $T \in \SD(X)$.
\item \label{r5}
If $T \in \Nar(X,\Gamma)$, then $T^{\cal U} \in \rNar(X,\Gamma^{\cal U})$.
\item \label{r6}
If $T^{\cal U} \in \Nar(X,\Gamma^{\cal U})$, then $T \in \Nar(X,\Gamma)$.
\end{enumerate}
\end{lemma}

\begin{proof}
All these statements don't differ too much in essence. Let us prove
for example \kref{r5}.  Fix arbitrary 
$x = (x_n)$, $ y  = (y_n) \in S(X^{\cal U})$ and 
$x^* = (x^*_n) \in \Gamma^{\cal U}$. Without loss of generality
(just replacing one representation of an element in $X^{\cal U}$ 
by another) one may assume that  $x_n, y_n \in S(X)$ for all 
$n \in \N$. Applying the condition $T \in \Nar(X,\Gamma)$ for 
 $x_n, y_n, x^*_n$ and $\eps = \frac{1}{n}$ we obtain elements 
$z_n \in S(X) $ such that $\|x_n+z_n\|> 2 - \frac{1}{n}$ and
$\|T(y_n-z_n)\|+|x^*_n(y_n-z_n)|< \frac{1}{n}$. This means that for
$z  = (z_n) \in S(X^{\cal U})$ the conditions $\|x+z\|= 2$ and
$\|T^{\cal U}(y-z)\|+|x^*(y-z)|=0$ are fulfilled.
\end{proof}


\section{Strong Daugavet and narrow operators in $\ell_\infty$-sums}
\label{sec3}

We first fix some notation. If $T$ is an operator defined on
$X=X_{1}\oplus_{\infty} X_{2}$, we let $T_{1}$ stand for the
restriction of $T$ to $X_{1}$, i.e., $T_{1}x_{1}= T(x_{1},0)$; and
likewise $T_{2}x_{2}= T(0,x_{2})$ defines the restriction to $X_{2}$. 
Thus for $x=(x_1,x_2) \in X$,
$Tx=T(x_1,x_2)=T_1x_1+T_2x_2$.

The aim of this section is to prove that $T$ is a
strong Daugavet operator if and only if both restrictions
$T_1$ and $T_2$ of $T$ are strong Daugavet operators. 
The same is true for narrow operators.

\begin{prop}\label{prop4.1}
If $X=X_1 \oplus_\infty X_2$ and $T_i \in \SD(X_i)$ 
($T_i \in \rSD(X_i)$) for $i=1,2$,
then $T \in \SD(X)$ ($T \in \rSD(X)$ respectively).
\end{prop}

\begin{proof}
 By Lemma~\ref{rigid} it is sufficient to consider only
the ``rigid'' version of the proposition. Indeed, we have $X^\calU =
X_{1}^\calU \oplus_{\infty} X_{2}^\calU$ and $(T^\calU)_{i}=
(T_{i})^\calU$. Therefore, if $T_{i}\in \SD(X_{i})$, then
$(T_{i})^\calU \in \rSD(X_{i}^\calU)$ and, assuming the rigid
version, we conclude that $T^\calU \in \rSD(X^\calU)$ which implies
$T\in \SD(X)$.

Thus, we need to prove that for every $x=(x_1,x_2)$ with 
$\|x\|=\max\{\|x_1\|,\allowbreak \|x_2\|\}=1$ and $y=(y_1,y_2)$ with 
$\|y\|=\max\{\|y_1\|,\|y_2\|\}=1$, there is some $z=(z_1,z_2)$ with 
$\|z\|=\max\{\|z_1\|,\|z_2\|\}=1$ such that
$\|x+z\|=\max\{\|x_1+z_1\|,\allowbreak \|x_2+z_2\|\}=2$ and $\|T(y-z)\| =
\|T_1(y_1-z_1)+T_2(y_2-z_2)\|=0$.

Without any loss of generality we may assume that $\|x_1\|=1$.
Using Lemma~\ref{equiv} for
$T_1 \in \rSD(X_1)$, we can find,
given $x_1 \in S(X)$ and $y_1 \in B(X)$, 
some $z_1 \in S(X)$ with  $\|x_1+z_1\|=2$ and
$\|T_1(y_1-z_1)\|=0$.
Put $z_2=y_2$, $z=(z_1,z_2)$; then $ \|z\|=1$ and
$\|x+z\| \geq \|x_1+z_1\|=2$ and
$$
\|T(y-z)\|=\|T_1(y_1-z_1)+T_2(y_2-z_2)\|=\|T_1(y_1-z_1)\|=0,
$$
completing the proof.
\end{proof}

\begin{cor} \label{n6}
If $X=X_1 \oplus_\infty X_2$  and $T_i \in
\Nar(X_i)$ for $i=1,2$, then $T \in \Nar(X)$.
\end{cor}

\begin{proof}
We have to prove that for each $x^*=(x_1^*,x_2^*) \in X^*=
X_{1}^* \oplus_{1} X_{2}^*$,
$T \mytilde{+} x^*$ is a strong Daugavet operator; see
Section~\ref{sec2} for this notation.

Let us consider the restriction of $T \mytilde{+} x^*$ to $X_1$;
then 
\begin{align*}
\|(T \mytilde{+} x^*)_1x_1\| 
&= \|(T \mytilde{+} x^*)(x_1,0)\| \\
&=
\|T(x_1,0)\|+|x^*((x_1,0))| \\
&=
\|T_1x_1\|+|x^*_1(x_1)|.
\end{align*}

Since $T_1$ is narrow, $T_{1} \mytilde+ x_{1}^*$ is a \sDo\ and hence
so is $(T\mytilde+ x^*)_{1}$.
By symmetry, the same is true for
the restriction to $X_{2}$, and
Proposition~\ref{prop4.1} implies
that $T \mytilde{+} x^*$ is a strong
Daugavet operator. Since $x^*$ is arbitrary, $T$ is narrow.
\end{proof}

We now turn to the converse of Proposition~\ref{prop4.1}.
In the sequel we shall call elements $x_1,\dots ,x_n$ of a normed
space {\em quasi-collinear}\/ if 
$$
\|x_1+\dots +x_n\|=\|x_1\|+\dots +\|x_n\|.
$$

We will need a simple lemma.

\begin{lemma}\label{lemma4.3}
Suppose that $x_1,\dots ,x_n$ are quasi-collinear.
\begsta
\item
$\|a_1x_1+\dots +a_nx_n\|=a_1\|x_1\|+\dots +a_n\|x_n\|$ for all
nonnegative coefficients $a_k$.
\item
If $x_{n+1}$ is quasi-collinear to $(x_1+\dots +x_n)/n$, then all the
vectors $x_1,\dots ,x_{n+1}$  are quasi-collinear.
\endsta
\end{lemma}

\begin{proof}
(a) The function $F\dopu \R_{+}^n\to \R$ defined by
$$
F(a_{1},\dots, a_{n})=
\|a_1x_1+\dots +a_nx_n\| - (a_1\|x_1\|+\dots +a_n\|x_n\|)
$$
is convex, takes values ${\le0}$ and $F(1,\dots,1)=0$. Hence $F=0$.

(b) follows from (a):
\bea
\|x_{1} + \dots + x_{n} + x_{n+1} \| &=&
\Bigl\| n \frac{x_{1} + \dots + x_{n} }{n}      + x_{n+1} \Bigr\| \\
&=&
n \Bigl\|  \frac{x_{1} + \dots + x_{n} }{n} \Bigr\| + \|x_{n+1} \| \\
&=&
\|x_{1} \|+ \dots + \|x_{n}\|   + \|x_{n+1} \| .
\eea
\end{proof}

\begin{thm} \label{inftysum}         
If $X= X_{1}\oplus_{\infty} X_{2}$, then for every strong 
Daugavet operator $T$ on $X$
the restrictions $T_1$ and $T_{2}$ of $T$ to $X_1$ and $X_{2}$
are strong Daugavet  operators.
\end{thm}

\begin{proof}
As in Proposition~\ref{prop4.1}
it is sufficient to prove that $T_1\in \SD(X_1)$ whenever
 $T \in \rSD(X)$.

So let $T \in \rSD(X)$, $x_1, y_1 \in S(X_1)$ and  $\eps > 0$. 
Apply the definition of a rigid strong Daugavet operator to
$x=(x_1,0)$, $y=(y_1,0)$. We get some $z^1=(z_1^1,z_2^1)$ for which 
$\|y_1+z_1^1\|=1$, $\|z_2^1\| \le 1$, $\|x_1+y_1+z_1^1\|=2$ and $Tz^1=0$. 
This means, in particular, that the vectors $x_1$ and $y_1+z_1^1$ 
are quasi-collinear. Now apply the definition of  a rigid
strong Daugavet operator to
$x=((x_1+y_1+z_1^1)/2,0)$, $y=(y_1,z_2^1)$. We get some $z^2=(z_1^2,z_2^2)$ 
for which $Tz^2=0$, $\|y_1+z_1^2\|=1$, $\|z_2^1+z_2^2\| \le 1$ and 
$\|(x_1+y_1+z_1^1)/2+(y_1+z_1^2)\|=2$. This again means, by
Lemma~\ref{lemma4.3}, that the 
vectors $x_1$, $y_1+z_1^1$  and $y_1+z_1^2$  are quasi-collinear.
Now apply the same token  to $x=((x_1+(y_1+z_1^1)+(y_1+z_1^2))/3,0)$ and 
$y=(y_1,z_2^1+z_2^2)$, etc.

Continuing this process we obtain a sequence $z^n=(z_1^n,z_2^n)$ 
for which all the vectors $x_1, y_{1}+ z_{1}^1,  y_1+z_1^2, \dots$
are quasi-collinear unit vectors, $\|z_2^1+\dots +z_2^n\| \le 1$  and $Tz^n=0$.
Consider $z=(z_1^1+z_1^2+\dots +z_1^n)/n \in X_1$. By construction
and Lemma~\ref{lemma4.3}
$\|x_1+y_1+z\|=2$, $\|y_1+z\|=1$ and 
$$
\|T_1z\|=\|T(z,0)\|=\|T(0,(z_2^1+z_2^2+\dots +z_2^n)/n)\| \le
\|T\|/n. 
$$
Because $n$ can be taken arbitrarily big, this proves that $T_1\in \SD(X_1)$.
\end{proof}

\begin{cor}\label{cor4.5}
If $X= X_{1}\oplus_{\infty} X_{2}$, then for every narrow
operator $T$ on $X$
the restrictions $T_1$ and $T_{2}$ of $T$ to $X_1$ and $X_{2}$
are narrow operators.
\end{cor}

\begin{proof}
This follows directly from Theorem~\ref{inftysum} and the definition
of a narrow operator.
\end{proof}

Let $X_1$ be an $M$-ideal of a Banach space $X$ (see \cite{HWW} for a
study of $M$-ideals) and $T$ be a
strong Daugavet operator on $X$. 
We haven't been able to decide whether 
 the restriction of
$T$ to $X_1$ is a strong Daugavet operator again. 
This would give us the operator version of the result  saying that an
$M$-ideal in a space with the Daugavet property has the 
Daugavet property itself \cite[Prop~2.10]{KadSSW}.


\section{Strong Daugavet and narrow operators in
  $\ell_1$-sums}
\label{sec4}

We use the same notation concerning restrictions of operators as 
before, but for an $\ell_{1}$-sum
$X=X_{1}\oplus_{1} X_{2}$.

\begin{prop}\label{prop5.1}
If $X=X_1 \oplus_1 X_2$ and $T_i \in \SD(X_i)$ ($T_i \in \rSD(X_i)$)
for $i=1,2$, then $T \in \SD(X)$ ($T \in \rSD(X)$ respectively).
\end{prop}

\begin{proof}
 Again, by Lemma~\ref{rigid} it is sufficient to consider only
the ``rigid'' version of the theorem.
Thus, we need to prove that
 for every $x=(x_1,x_2)$ with 
$\|x\|=\|x_1\|+\|x_2\|=1$ and $y=(y_1,y_2)$ with 
$\|y\|=\|y_1\|+\|y_2\|=1$, there is some $z=(z_1,z_2)$ with
$\|z\|=\|z_1\|+\|z_2\|=1$ such that
$\|x+z\|=\|x_1+z_1\|+\|x_2+z_2\|=2$ and $\|T(y-z)\| =
\|T_1(y_1-z_1)+T_2(y_2-z_2)\|=0$.

For $i=1,2$, since $T_i \in \rSD(X_i)$, 
we can produce, using Lemma~\ref{equiv}, some $z_i \in \|y_i\|S(X_i)$ with
$\|x_i+z_i\|=\|x_i\|+\|z_i\|$ and
$\|T_i(y_i-z_i)\|=0$.
Now let us take $z=(z_1,z_2)$; then
\begin{eqnarray*}
&\|z\|=\|z_1\|+\|z_2\|=\|y_1\|+\|y_2\|=1,&
\\
 &\|x+z\|=\|x_1+z_1\|+\|x_2+z_2\|=\|x_1\|+\|z_1\|+\|x_2\|+\|z_2\|
=2&
\end{eqnarray*}
and
$$
T(y-z)=T_1(y_1-z_1)+T_2(y_2-z_2)=0.
$$
 So, $z$ satisfies all the conditions above, and the proposition is
proved.
\end{proof}

By the same argument as in Corollary~\ref{n6}
we obtain:

\begin{cor}\label{cor5.2}
If $X=X_1 \oplus_1 X_2$  and $T_i \in
\Nar(X_i)$ for  $i=1,2$, then $T \in \Nar(X)$.
\end{cor}

We now study the converse of these results.
Let us recall that a subset $\Gamma \subset S(X^*)$ is said to be 
$1$-norming if
$$
\|x\|= \sup _{x^*\in \Gamma} x^*(x).
$$
 for every $x \in X$.
A subset $\Gamma \subset S(X^*)$ is said to be a boundary for $X$ 
if the above supremum is always attained, i.e.,
if for every  $x \in X$ there is some $x^* \in \Gamma$ such that  
$x^*(x) = \|x\|$. Clearly, the notion of a boundary is 
a ``rigid'' version of a $1$-norming set. It is easy 
to check that  $\Gamma^{\cal U}$ is a boundary
for  $X^{\cal U}$ if and only if  $\Gamma$ is $1$-norming.

\begin{lemma}
Let $X=X_1 \oplus_1 X_2$, let $\Gamma_j \subset S(X^*_j)$ be
boundaries for $X_j$ for
$j=1,2$, and let $\Gamma =\Gamma_1 \cup \Gamma_2 $. 
If  $T \in \rNar(X,\Gamma)$, then
$T_j$, the restrictions of $T$ to $X_j$, are rigid  strong
Daugavet operators.
\end{lemma}

\begin{proof}
Let us consider the case of $T_1$.
We have to prove that for every  
$x_1,y_1 \in S(X_1)$ there exists some $u_1 \in S(X_1)$ such that
$\|x_1+u_1\|=2$ and $T_1(u_1-y_1)=0$.

Let us take $x=(x_1,0)$, $ y=(y_1,0) \in S(X)$ and a 
functional $x^*_1 \in \Gamma_1$
such that $x^*_1(y_1)=1$. Let us further take $x^*=(x_1^*,0) \in \Gamma$.
Since $T$ is narrow, we can apply Definition~\ref{rigiddef}  with
the elements $x,y$ and $x^*$ defined above; thus, 
there exists some $z=(z_1, z_2) \in S(X)$
such that
$$\|x+z\|=\|x_1+z_1\|+\|z_2\|=2$$
and
\begin{equation} \label{66}
\|T(z-y)\|+|x^*(z-y)|=\|T(z-y)\|+|x^*_1(z_1-y_1)|=0.
\end{equation}
 From the last condition we obtain
$|x^*_1(z_1-y_1)|=0$. Keeping in mind that
$x^*_1(y_1)=1$, we get $x^*_1(z_1)=1$. But
$\|x^*_1\|=1$, so $\|z_1\|=1$. Then
\begin{equation} \label{666}
\|z_2\|=0, 
\end{equation}
because $\|z_1\|+\|z_2\|=1$.
So
$\|x_1+z_1\| = \|x+z\|=2$
and by (\ref{66}) and (\ref{666})
$T_1(y_1)=T(y)=T(z)=T_1(z_1)$.
Thus the definition of a rigid strong Daugavet operator 
is fulfilled for $T_1$.
\end{proof}

We can now prove the converse of Corollary~\ref{cor5.2}.

\begin{thm}\label{theo5.4}
Let $X=X_1 \oplus_1 X_2$ and $T \in \Nar(X)$. Then $T_1$ and $T_{2}$, the
restrictions of $T$ to $X_1$ and $X_{2}$, are
narrow operators.
\end{thm}

\begin{proof}
 It has been proved in \cite[Cor.~3.14]{KadSW2} that if $T$ is
 narrow then so is $T \mytilde{+} x^*$ for any 
$x^* \in  X^* $, in particular for $x^*\in \Gamma = X_{1}^* \cup
X_{2}^*$. By Lemma~\ref{rigid} we may pass to
ultraproducts, apply  the previous lemma, pass back to the original
space and obtain that
$T_1 \mytilde{+} x_1^*$ is strongly Daugavet for every $x_1^* \in  X_1^*$.
Hence, by definition, $T_1$ is narrow, and by symmetry, so is $T_{2}$.
\end{proof}

However, the analogue of Theorem~\ref{theo5.4} for \sDo s, i.e., the
converse of Proposition~\ref{prop5.1}, is false.

\begin{prop}
Let $X=X_1 \oplus_1 X_2$ and $T \in \SD(X)$. Then $T_1$, the
restriction of $T$ to $X_1$, need not be a \sDo.
\end{prop}

\begin{proof}
The sum functional $Tx= \sum_{n=1}^\infty x(n)$ is a \sDo\ on
$\ell_{1}= \R \oplus_{1} X_{2}$ (see \cite[Prop.~2.4]{BKSSW}), 
yet its restriction to
$\R$ (i.e., the span of $e_{1}$) is not.
\end{proof}

We wish to indicate another counterexample that even works on a space
with the \DP, viz.\ $L_{1}[0,1]$. For this,
let us recall the main features of the example from 
Theorem~6.3 of \cite{KadSW2}. In this example  subspaces 
$Y_1 \subset L_1[0,1]$ and $Y= Y_{1} \oplus \lin\{ 1\}$ 
and a measurable subset 
$P \subset [0,1]$ of measure
$\mu(P) < 1/9 $  with the following properties are constructed:
\begin{equation} \label{uh}
\|g \chi_{[0,1] \setminus  P} \| \le 3 \|g \chi_{P}\| \qquad \forall
g \in Y_1
\end{equation}
and
the quotient map $q\dopu   L_1[0,1] \rightarrow L_1[0,1]/Y$ is a  
strong Daugavet operator.

Now let $Q \subset [0,1]$, $\mu(Q) < 1/3$,
$Q \cap P = \emptyset$. Then the restriction of $q$
to $ L_1(Q)$ is bounded from below. So in particular this 
restriction is not a strong Daugavet operator; observe that
$L_{1}[0,1] = L_{1}(Q) \oplus_{1} L_{1}([0,1]\setminus Q)$.

Indeed, let us assume to the contrary that the restriction of $q$
to $ L_1(Q)$ is unbounded from below. This means that for every
$\eps > 0$ there exists a function $f \in  L_1(Q)$, a function
$g_1 \in Y_1$ and a constant $a$ such that 
$$
\|f - (g_1+a)\| < \eps.
$$
Denote $[0,1] \setminus ( P \cup Q)$ by $S$; then $\mu(S) >
1/2$. 
Then $ \|(a+g_1) \chi_{P \cup S} \| < \eps$ and
\begin{align*}
a \mu(P) = 
\|a \chi_{P} \| &\ge  
\|g_1 \chi_{P}\| - \eps 
\ge 
\frac{1}{3}\|g_1 \chi_{S}\| - \eps 
\qquad\text{(by~\ref{uh})}\\
&\ge  
\frac{1}{3}\|a \chi_{S} \| - 2 \eps = 
\frac{1}{3} a \mu(S) - 2 \eps,
\end{align*}
so $a < 40 \eps$.
This means that $\|f - g_1\| < 41\eps$. On the other hand
\begin{align*}
\|f - g_1\| 
&\ge 
\|(f-g_1) \chi_{P}\| = 
\|g_1 \chi_{P}\| \\
&\ge 
\frac{1}{3}\|g_1 \chi_{Q}\| \ge 
\frac{1}{3}(\|f\| - \|(f-g_1) \chi_{Q}\|) \ge 
\frac{1}{3}(1 - 41 \eps),
\end{align*}
which is a contradiction.


\section{The Daugavet property for unconditional sums of spaces}
\label{sec5}

Throughout this section $F$ denotes a Banach space with a
$1$-unconditional normalised Schauder basis. We can think of the
elements of $F$ as sequences with the property that
$$
\| (a_{1}, a_{2}, \dots) \|_{F} =
\| (|a_{1}|, |a_{2}|, \dots) \|_{F} \qquad\forall (a_{j})\in F.
$$
Note that $F$ is naturally endowed with the structure of a Banach
lattice with respect to the pointwise operations. 

Suppose that $X_{1}, X_{2}, \dots$  are Banach spaces. Their $F$-sum
$X= (X_{1},\allowbreak X_{2}, \dots)_{F}$ consists of all sequences $(x_{j})$
with $x_{j}\in X_{j}$ and $(\|x_{j}\|)\in F$ with the norm
$\|(x_{j})\| = \|(\|x_{j}\|)\|_{F}$. We are going to characterise
when such an $F$-sum has the \DP. 

\begin{thm}\label{theo7.1}
Let $X_{1},X_{2},\dots$ be Banach spaces with the \DP. Then their
$F$-sum $X$ has the  \DP\ if and only if the Banach lattice $F$ has
the positive \DP\ in the sense that $\|\Id+T\|=1+\|T\|$ whenever
$T\dopu F\to F$ is a positive rank-$1$ operator. 
\end{thm}

\begin{proof}
We first remark that the positive \DP\ may be characterised as in
Lemma~\ref{l:MAIN}; the proof is verbatim the same as in
\cite[Lemma~2.1]{KadSSW}.

\begin{lemma}\label{lemma7.2}
A Banach lattice has the positive \DP\ if and only if for every
positive $a\in S(F)$, every positive $a^*\in S(F^*)$ and every
$\eps>0$ there is some positive $b\in S(F)$ such that $a^*(b)\ge 1-\eps$  and
$\|a+b\| \ge 2-\eps$.
\end{lemma}

Now suppose that $X$ has the \DP; we shall verify the condition of
Lemma~\ref{lemma7.2}. Note that $F^*$ can be represented by all
sequences $(a_{j}^*)$ such that 
$$
\sup_{n} \bigl\| ( |a_{1}^*| ,  \dots , |a_{n}^*| ,0,0,\dots) \bigr\|_{F^*}
<\infty,
$$
and $X^*$ can be represented by all
sequences $(x_{j}^*)$, $x_{j}^*\in X_{j}^*$,  such that 
$$
\|x^*\| =
\sup_{n} \bigl\| ( \|x_{1}^*\| ,  \dots , \|x_{n}^*\| ,0,0,\dots) 
\bigr\|_{F^*} <\infty.
$$
Let $a=(a_{j})\in S(F)$ and $a^*= (a_{j}^*) \in S(F^*)$ be positive
elements and let $\eps>0$. Pick $x_{j}\in X_{j}$ and $x_{j}^* \in
X_{j}^*$ such that $\|x_{j}\|=a_{j}$, $\|x_{j}^*\|= a_{j}^*$ and put
$x=(x_{j})$, $x^*= (x_{j}^*)$; then $\|x\| = \|x^*\|=1$. Since $X$
has the \DP, we can find $y\in S(X)$ such that $x^*(y)\ge 1-\eps$ and
$\|x+y\|\ge 2-\eps$; cf.\ Lemma~\ref{l:MAIN}. Write $y=(y_{j})$ and
$b=(\|y_{j}\|)$; then $\|b\|_{F}=1$ and
\bea
&\displaystyle
1-\eps \le x^*(y) = \sum_{j=1}^\infty x_{j}^*(y_{j}) \le \sum_{j=1}^\infty
\|x_{j}^*\| \|y_{j}\| = a^*(b), &\\
&2-\eps \le \|x+y\| =
\bigl\| (\|x_{j} + y_{j}\| ) \bigr\|_{F} \le
\bigl\| (\|x_{j} \| + \| y_{j}\| ) \bigr\| _{F} \le
\|a\| + \|b\| , &
\eea
where we have used the fact that the norm of $F$ is monotonic in each
variable. Hence $F$ has the positive \DP. (Incidentally, the
assumption that the $X_{j}$ have the \DP\ did not enter this part of
the proof.)

Conversely, suppose that $F$ has the positive \DP. Let $x=(x_{j}) \in
S(X)$  and $x^*= (x_{j}^*)\in S(X^*)$, define $a= (a_{j}) =
(\|x_{j}\|) \in S(F)$ and $a^*= (a_{j}^*) =
(\|x_{j}^*\|) \in S(F^*)$. Given $\eps>0$, find using
Lemma~\ref{lemma7.2} some $b=(b_{j})\in S(F)$ such that $a^*(b) \ge
1-\eps$ and $\|a+b\| \ge 2-\eps$. Since $X_{j}$ has the \DP, one can
find $y_{j}\in X_{j}$ such that $\|y_{j}\|= b_{j}$,
$x_{j}^*(y_{j}) \ge (1-\eps) a_{j}^* b_{j}$ and
$\|x_{j} + y_{j}\| \ge (1-\eps) (a_{j}+b_{j})$;
just note that $\|\Id + (x_{j}^*/a_{j}^*) \otimes (x_{j}/b_{j})\| =
1+a_{j}/b_{j}$. Therefore $y=(y_{j}) \in S(X)$ satisfies
$$
x^*(y) = \sum_{j=1}^\infty x_{j}^*(y_{j}) \ge
(1-\eps) \sum_{j=1}^\infty a_{j}^*b_{j} =
(1-\eps) a^*(b) \ge (1-\eps)^2
$$
and 
\bea
\|x+y\| &=&
\bigl\| (\|x_{j}+ y_{j}\|) \bigr\|_{F} 
\ge
(1-\eps) \bigl\| (\|x_{j}\| + \|y_{j}\|) \bigr\|_{F} \\
&=&
(1-\eps) \|a+b\|_{F} \ge 2 (1-\eps)(1-2\eps).
\eea
Hence $X$ has the \DP.
\end{proof}

It is clear that for example $c_{0}$ and $\ell_{1}$ have the positive
\DP, hence Theorem~\ref{theo7.1} contains \cite[Prop.~2.16]{KadSSW} as
a special case. Also, if $F^*$ is a Banach lattice with the positive
\DP, then so is $F$. 

If $F$ is finite-dimensional, we can pass to the limit $\eps=0$ in
Lemma~\ref{lemma7.2} by compactness. Thus, we obtain the following
variant of Theorem~\ref{theo7.1}.

\begin{cor} \label{cor7.3}
Let $\dim F=n$ and $X_{1},\dots, X_{n}$ be Banach spaces with the
\DP. Then their $F$-sum $(X_{1}\oplus \dots \oplus X_{n})_{F}$ has
the \DP\ if and only if for every positive $a\in S(F)$ and every
positive $a^*\in S(F^*)$ there is some $b\in S(F)$ such that
$a^*(b)=1$ and $\|a+b\|=2$.
\end{cor}

This condition can be rephrased geometrically as follows. For any
point $a\ge 0$ in $S(F)$ and any supporting hyperplane  of
the positive part of the
unit sphere $H=\{ a^*=1 \}$ there is a line segment in the unit
sphere that contains $a$ and intersects $H\cap S(F)$. From this the
following corollary is evident.

\begin{cor} \label{cor7.4}
If $X= (X_{1} \oplus X_{2})_{F}$ has the \DP, then either
$F=\ell_{1}^2$ or $F=\ell_{\infty}^2$, i.e., either 
$X= X_{1} \oplus_{ 1} X_{2}$ or $X= X_{1} \oplus_{ \infty} X_{2}$.
\end{cor}

It is easy to see that $F_{1} \oplus _{1 } F_{2}$ and
$F_{1} \oplus _{\infty } F_{2}$ have the positive \DP\ whenever $F_{1}$
and $F_{2}$ have; in fact, the proof of Theorem~\ref{theo7.1} shows
that the $F$-sum $(F_{1} \oplus F_{2} \oplus \dots)_{F}$ of Banach
lattices with the positive \DP\ is a Banach lattice with the positive
\DP. Therefore, starting from the real line we can form
$\ell_{1}$-sums and $\ell_{\infty}$-sums consecutively to obtain
finite-dimensional spaces with the positive \DP, e.g., the
$18$-dimensional space
$$
(\ell_{\infty}^3 \oplus_{1} \ell_{\infty}^4) \oplus_{\infty} 
(\ell_{1}^3 \oplus_{1} \ell_{\infty}^3) \oplus_{\infty} \ell_{1}^5.
$$
However, there are other examples, even in the three-dimensional
case; for
$$
\|(a_{1}, a_{2}, a_{3})\|_{F}  =
\max \Bigl\{ |a_{1}| + \frac{|a_{3}|}2 , |a_{2}|+|a_{3}| \Bigr\}
$$
defines a norm on $\R^3$ with the positive \DP. In this example the
unit sphere intersected with the half-space $\{(0,0,t)\dopu t\ge0\}$
looks like a hip roof and the positive part of $B(F)$, i.e.,
$B(F)\cap \R_{+}^3$, is the convex hull of the points
$(0,0,0)$, $(1,0,0)$, $(0,1,0)$, $(0,0,1)$, $(1,1,0)$ and
$(1/2,0,1)$. From this description it is easy to see (literally) that
this norm has the positive \DP.


\section{The Daugavet property for ultraproducts}
\label{sec6}

Let $Z$ be a subspace in $X^*$, $\Gamma = S(Z)$ be a boundary for $X$, and  
$X \in \DPr(\Gamma)$. Is it true that under this condition
$X$ has the \DP?
Provided the answer to this question is positive 
Lemma~\ref{rigid} easily implies that an ultrapower of a
space with the Daugavet property has the Daugavet property
itself. Unfortunately we don't know the answer;
that is why we investigate
the question from another point of view in this section.
Note, however, that it is easy to find a Banach space without the \DP\
that has the \DP\ with respect to some boundary, e.g., $\ell_1$ with
the boundary $\ex B_{\ell_1^*}$. 

For an element $x\in S(X)$ and $\eps>0$ denote 
$$
l^+(x,\eps) = \{y\in X\dopu \|y\|\le 1+\eps,\  \|x+y\|>2-\eps \}.
$$

The next lemma follows
directly from Lemma~\ref{l:MAIN}.

\begin{lemma} \label{lem1.2}
The following assertions are equivalent:
\begin{enumerate}
\item
$X$ has the Daugavet property.
\item
For every $x\in S(X)$ and every $\eps>0$ 
the closure of $\conv(l^+(x,\eps))$ contains
 $B(X)$.
\end{enumerate}
\end{lemma}

Lemma~\ref{lem1.2} suggests the following quantitative approach to
the \DP. For a subset $A \subset X$ denote by 
$\conv_n(A)$ the set of all convex 
combinations of all $n$-point collections of elements of $A$.
Clearly
$\conv(A)= \bigcup_{n \in \N} \conv_n(A)$. Denote  
$$ 
\Daug_n(X,\eps)=\sup_{x,y\in S(X)}\dist(y, \conv_n(l^+(x,\eps)))
$$
It is easy to see that for every $\eps>0$ the sequence
$\bigl( \Daug_n(X,\eps) \bigr)$ decreases. If 
$\lim_{n\to \infty} \Daug_n(X,\eps) =0$
for every $\eps>0$, 
then $X$ has the Daugavet property. 
We don't know if the converse is true.  

\begin{thm} \label{ultraproduct}
Let $\cal U$ be an ultrafilter defined on a set ${\cal I}$, and let $X_{i }$,
$i  \in {\cal I}$,
be a collection of Banach spaces and $X$ be the corresponding
ultraproduct of the $X_{i }$. Then the following assertions are
equivalent:
\begin{enumerate}
\item \label{A}
$X$ has the Daugavet property.
\item \label{B}
For every $\eps>0$, $\lim_{\cal U,n}\Daug_n(X_{i },\eps)=0$. In
other words, for every fixed $\eps>0$ and  every
$\delta>0$ there is  an $n \in \N$ such
that the set of all $i $ for which $\Daug_n(X_{i },\eps) < \delta$ 
belongs to the  ultrafilter $\cal U$.
\end{enumerate}
\end{thm}

\begin{proof}
To deduce~\kref{A}  from~\kref{B} one just has to notice that if the 
set of all $i $ for which $\Daug_n(X_{i },\eps) < \delta$ 
belongs to the  ultrafilter $\cal U$, then $\Daug_n(X,\eps) < \delta$.
So  for every $\eps>0$, $\Daug_n(X,\eps)$ tends
to 0 when $n$ tends to infinity, which proves the Daugavet property
for $X$.

To deduce~\kref{B} from~\kref{A} let us argue ad absurdum. Suppose there
are  $\eps>0$ and  $\delta>0$ such that for every $n \in \N$ the set
$A_n=\{i  \in {\cal I}\dopu  \Daug_n(X_{i },\eps) > \delta\}$ belongs 
to the  ultrafilter $\cal U$. Denote $A_0 = {\cal I}$. Let us construct
two elements $x=(x_{i })_{i  \in {\cal I}}$ and 
$y=(y_{i })_{i  \in {\cal I}}$ of $S(X)$ in such a way that 
$x_{i },y_{i } \in S(X_{i }))$ and for
every $i  \in A_n \setminus A_{n-1}$ the distance from 
$\conv_n(l^+(x_{i },\eps))$ to $y_{i }$ is bigger than $\delta$. 
The $\conv_n$-hull of a set is
increasing when $n$ is increasing, so for every $n \in \N$
and every $i  \in A_n= \bigcup_{m=n}^{\infty}A_m \setminus A_{m-1}$ 
the distance from $\conv_n(l^+(x_{i },\eps))$ to $y_{i }$ 
is bigger than $\delta$. This means in turn  that for every $n \in
\N$, 
$\dist(y, \conv_n(l^+(x,\eps))) \ge \delta$, so 
$\dist(y, \conv(l^+(x,\eps))) \ge \delta$, which contradicts the 
Daugavet property of $X$.
\end{proof}

\begin{rem}\label{rem6.3}
If 
$\lim_{n\to \infty} \Daug_n(X,\eps) =0$
for every $\eps>0$, 
then for every $\eps>0$
there is some $n \in \N$ such that $\Daug_n(X,\eps)=0$. 
More explicitly:
If  $\Daug_n(X,\eps/2)< \eps/2$, then $\Daug_n(X,\eps)=0$. Moreover
for every pair $x,y\in S(X)$ not just
$\dist(y, \conv_n(l^+(x,\eps)))=0$, but $y\in \conv_n(l^+(x,\eps))$.
\end{rem}

\begin{proof}
Suppose $\Daug_n(X,\eps/2)< \eps/2$. Fix $x,y\in S(X)$. There 
exist $y_1, \dots  , \allowbreak y_n \in (1+\eps/2)B(X)$,
$\|x+y_n\|>2-\eps/2$, and $a_1, \dots , a_n \ge 0, \sum_{k=1}^n a_k =1$, 
for which $\|y- \sum_{k=1}^n a_k y_k\| < \eps/2$.
Define elements $z_j= y_j + y - \sum_{k=1}^n a_k y_k$. Then  
$z_j \in l^+(x,\eps)$, $\sum_{j=1}^n a_j z_j = y$, so $y\in
\conv_n(l^+(x,\eps))$.
\end{proof}

So instead of $\Daug_n(X,\eps)$ it is reasonable to consider the
following notion, which seems to be a bit more convenient (at least
it depends only on one parameter):
$$
D_X(\eps)=\inf\{n\dopu  \conv_n(l^+(x,\eps)) 
\supset  S(X)\ \forall x \in S(X) \}
$$

If $D_X(\eps)$ is finite for every $\eps>0$, we say that $X$
possesses the \textit{uniform Daugavet property}.
Equivalently, by Remark~\ref{rem6.3}, $X$ has the uniform \DP\ if and
only if $\Daug_{n}(X,\eps)\to 0$ for every $\eps>0$.

The theorem on ultraproducts can be reformulated in the following
way.

\begin{thm} \label{ultraproduct11}
Let $\cal U$ be an ultrafilter defined on a set ${\cal I}$, $X_{i }$
be a collection of Banach spaces and $X$ be the corresponding
ultraproduct of the $X_{i }$. Then the following assertions are
equivalent:
\begin{enumerate}
\item \label{A11}
$X$ has the Daugavet property.
\item\label{B11}
For every $\eps>0$ there exists some $n$ such that the set of all 
$i $ for which $D_{X_{i }}(\eps)<n$
belongs to the  ultrafilter $\cal U$.
\end{enumerate}
\end{thm}

\begin{cor}
A Banach space $X$ has the uniform \DP\ if and only if  
every ultrapower $X^{\cal U}$ has the \DP, in which case 
$X^{\cal U}$ even has the uniform \DP.
\end{cor}

It follows from this corollary and the canonical isometric isomorphism
$(X \oplus_\infty  Y)^\calU = X^\calU \oplus_\infty
Y^\calU$ that the uniform \DP\ is stable by taking
$\ell_\infty$-direct sums and likewise by taking
$\ell_1$-direct sums.

Let us prove that the basic examples
of spaces with the Daugavet property
in fact are spaces with the uniform Daugavet property.

\begin{lemma}
Let $X=L_{1}[0,1]$. If $n > 2/\eps $, then 
$\Daug_{n}(X,\eps )=0$; if $n \le 2/\eps$, then 
$\Daug_{n}(X,\eps )\leq 1- \eps n/(2+\eps )$.
Hence $D_{X}(\eps)$ is of order $\eps^{-1}$.
\end{lemma}

\begin{proof}
Suppose $n > 2/\eps $ and let us take arbitrary points $x$ and $y$
from $S(X)$. There is a partitioning of $[0,1]$  into
sets $E_{1},\dots  ,E_{n}$ such that $\Vert x\cdot \chi _{E_{i}}\Vert =
1/n < \varepsilon /2$. Define functions $y_{i}$ by
$y_{i}=\frac{1}{\Vert y\cdot \chi _{E_{i}}\Vert } y\cdot \chi _{E_{i}}$ 
if $\Vert y\cdot
\chi _{E_{i}}\Vert \neq 0$, and $y_{i}=0$ if $\Vert y\cdot \chi
_{E_{i}}\Vert =0$. Then 
$\sum_{i=1}^{n}y_{i}\lambda _{i}=y$,
where $\lambda _{i}=\Vert y\cdot \chi _{E_{i}}\Vert $. On the other hand, if 
$y_{i}\neq 0$, then
$$
\Vert x+y_{i}\Vert  \geq \Vert x\cdot \chi _{[0,1]\backslash E_{i}}\Vert
+\Vert y_{i}\Vert -\Vert x\cdot \chi _{E_{i}}\Vert \geq 2-2\Vert x\cdot \chi
_{E_{i}}\Vert  > 2-\eps .
$$
So, $y_{i}\in l^{+}(x,\eps )$.

If $n \le 2/\eps $, then proceeding as above, with $N=[2/\eps ]+1$ 
we get a decomposition $E_{1},\dots , E_{N}$. Let us
arrange the $\lambda _{i}$'s in decreasing order 
and take the first $n$ of them.
Then 
$$
\biggl\Vert \sum_{i=1}^{n}y_{i}\lambda _{i}-y \biggr\Vert  =
\biggl\Vert
\sum_{i=n+1}^{N}y_{i}\lambda _{i} \biggr\Vert  
\leq \sum_{i=n+1}^{N}\lambda _{i}=S.
$$
We need to prove that $S\leq (N-n)/N$. Assume the opposite. Then 
$$
1 = \sum_{i=1}^{N}\lambda _{i}>\sum_{i=1}^{n}\lambda _{i}+\frac{N-n}{N}; 
$$
hence
$n/N > \sum_{i=1}^{n}\lambda _{i}\geq n\lambda _{n}$ and 
$1/N > \lambda _{n}$.
Thus, 
\[
S=\sum_{i=n+1}^{N}\lambda _{i}\leq \lambda _{n}(N-n)<\frac{N-n}{N}\text{,}
\]
which is a contradiction. 

So, 
$$
S\leq \frac{N-n}{N}=1-\frac{n}{[\frac{2}{\eps }]+1}\leq 1-\frac{%
\eps n}{2+\eps }
$$ 
 and the proof of the lemma is finished.
\end{proof}

\begin{lemma}
If $X=C(K)$ for a compact Hausdorff space $K$ without isolated
points, then for every $\eps $ and $n$, $\Daug_{n}(X,\eps
)\leq {2}/{n}$.
Hence $D_{X}(\eps)$ is of order $\eps^{-1}$.
\end{lemma}

\begin{proof}
Let $x$ and $y\in S(X)$ be arbitrary. Without loss of generality,
assume that $x$ attains the value~$1$.
Take an open neighbourhood $U$ such that $x(u)>1-\eps $ for all $u\in U
$. Now pick $n$ disjoint subneighbourhoods $V_{1},\dots  ,V_{n}$ inside $U$.
For each of them choose a positive function $\varphi _{i}$ supported on 
$V_{i}$ such that $\Vert \varphi _{i}\Vert \leq 2$, $\Vert y+\varphi
_{i}\Vert \leq 1$ and $y+\varphi _{i}$ attains the value~$1$ 
in $V_{i}$. Obviously, 
$\Vert x+y+\varphi _{i}\Vert >2-\eps $, hence, $y+\varphi _{i}\in
l^{+}(x,\eps )$. On the other hand,
\[
\biggl\| \frac{1}{n}\sum_{i=1}^{n}(y+\varphi _{i})-y \biggr\| 
= \biggl\| \frac{1%
}{n}\sum_{i=1}^{n}\varphi _{i}\biggr\| \leq \frac{2}{n}\text{,}
\]
which proves the lemma. 
\end{proof}

One can show that
the same estimates 
for the $\Daug_{n}$ constants are valid   for rich subspaces
of  $C(K)$-spaces (see \cite{KadSW2} for this notion), 
for vector-valued $C(K)$- or $L_1$-spaces
and for spaces of weakly continuous vector-valued functions with
the sup-norm.

\begin{rem}
For every Banach space $X$ the constants $D_{X}(\eps)$
can be estimated from below by 
$(2+2 \eps)/(3 \eps)$, which is bigger than $2/(3 \eps)$. 
So the estimates from 
above which we have for $L_1$ and $C$ are of optimal order.
\end{rem}

\begin{proof}
Suppose $D_{X}(\eps)=n< (2+2 \eps)/(3 \eps)$ for some $\eps>0$.
This means in particular that for a fixed element $x\in S(X)$ (taking 
$y= -x$) there
are elements $y_1, \dots  , y_n \in (1+\eps)B(X)$,
$\|x+y_n\|>2-\eps$ and $a_1, \dots  ,a_n \ge 0$, $ \sum_{k=1}^n a_k =1$, 
for which $ \sum_{k=1}^n a_k y_k=-x$. Without loss of generality we may
assume that $a_1 \ge 1/n$ (otherwise just change the enumeration).
Plugging in  $\|x+y_1\|>2-\eps$  and $ x= -\sum_{k=1}^n a_k y_k$ we obtain
\begin{align*}
2-\eps < \|x+y_1\| &= 
\biggl\|y_1(1-a_1) - \sum_{k=2}^n a_k y_k \biggr\| \\
&\le
(1+\eps)(1-a_1) + (1+\eps)(1-a_1) \\
&\le 
2(1+\eps)(1-1/n) \le 2-\eps,
\end{align*}
which is a contradiction.
\end{proof}



\end{document}